\documentclass[11pt]{amsart}

\usepackage{amsfonts,amsmath,amsthm,amssymb}
\usepackage{latexsym}

\newtheorem{theorem}{Theorem}[section]

\newtheorem{proposition}[theorem]{Proposition}

\theoremstyle{definition}

\newcommand{\Z}{\ensuremath{\mathbb{Z}}}
\newcommand{\R}{\ensuremath{\mathbb{R}}}

\newcommand{\Co}{\ensuremath{\mathbb{C}}}

\newcommand{\T}{\ensuremath{\mathcal{T}}}

\def \W {\Omega}

\def \< {\langle}
\def \> {\rangle}

\newcommand{\beql}[1]{\begin{equation}\label{#1}}
\newcommand{\eeq}{\end{equation}}
\newcommand{\comment}[1]{}

\newcommand{\Abs}[1]{{\left|{#1}\right|}}

\newcommand{\Set}[1]{{\left\{{#1}\right\}}}

\newcounter{rem}
\begin{document}

\title{Complex Hadamard matrices and the Spectral Set Conjecture}

\author{Mihail N. Kolountzakis \& M\'at\'e Matolcsi}

\date{October 2004}

\address{M.K.: School of Mathematics, Georgia Institute of Technology, 686 Cherry St NW, Atlanta, GA 30332, USA, and
Department of Mathematics, University of Crete, Knossos Ave.,
GR-714 09, Iraklio, Greece} \email{kolount@member.ams.org}

\address{M.M.: Alfr\'ed R\'enyi Institute of Mathematics,
Hungarian Academy of Sciences POB 127 H-1364 Budapest, Hungary.}\email{matomate@renyi.hu}

\thanks{Partially supported by European Commission IHP Network HARP
(Harmonic Analysis and Related Problems), Contract Number: HPRN-CT-2001-00273 - HARP}

\begin{abstract}
By analyzing the connection between complex Hadamard matrices and spectral
sets we prove the direction ``spectral $\Rightarrow$ tile'' of the Sectral Set
Conjecture for all sets $A$ of size $|A|\le 5$ in any finite Abelian group. This
result is then extended to the infinite grid $\Z^d$ for any dimension $d$,
and finally to $\R^d$.

It was pointed out recently in \cite{tao} that the
corresponding statement fails for $|A|=6$ in the group $\Z_3^5$, and this
observation quickly led to the failure of the Spectral Set Conjecture in
$\R^5$ \cite{tao}, and subsequently in $\R^4$ \cite{mat}.
In the second part of this note we reduce this dimension further, showing that the direction ``spectral
$\Rightarrow$ tile'' of the Spectral Set Conjecture is false already in dimension 3.

In a computational search for counterexamples in lower dimension (one and two) one needs,
at the very least, to be able to decide efficiently if a set is a tile (in, say, a cyclic group) and if
it is spectral. Such efficient procedures are lacking however and we make a few comments for
the computational complexity of some related problems.
\end{abstract}

\maketitle



\section{Introduction}\label{sec:intro}
Let $G$ be a locally compact Abelian group and $\Omega\subseteq G$ be a bounded open set.
We call $\Omega$ {\em spectral} if there exists a set $\Lambda$ of continuous characters of
$G$ which forms an orthogonal basis of the space $L^2(\Omega)$ (with respect to Haar measure).
Such a set $\Lambda$ is called a {\em spectrum} of $\W$.

The set $\Omega$ tiles $G$ by translation
if there exists a set $T\subseteq G$ of translates, called a {\em tiling complement} of $\Omega$,
such that $\sum_{t\in T} \chi_\Omega(x-t) = 1$, for almost all (with
respect to the Haar measure) $x\in G$.
Here $\chi_\Omega$ denotes the indicator function of $\Omega$.

The {\it Spectral Set Conjecture} (originally proposed in $G=\R^d$ by Fuglede,
\cite{fug}) states that $\W$ is a tile in $G$ if and only if it is spectral.

This conjecture has attracted considerable attention over the last
decade. Many positive results were obtained in special cases (i.e., under various restrictions on $\W,
G$, and $T$; see
\cite{fug,convexcurv,convexplane,spectralsym,konyagin,lagwang,lagszab,unispectra}),
until  Tao \cite{tao}  recently proved that the direction ``spectral $\Rightarrow$ tile'' does not
hold in dimension 5 and higher. Tao's result was supplemented by a result
of Matolcsi \cite{mat},  which allowed the dimension of the counterexample to
be reduced to 4. In the second part of this paper we combine the ideas of \cite{tao}
and \cite{mat} in order to reduce the dimension to 3.
We also remark that the ``tile $\Rightarrow$ spectral'' direction was recently disproved
in \cite{kolmat} in dimension 5 and higher.
The existing counterexamples do not exclude the possibility that the
conjecture may still be true in low dimensions (in particular, in dimension 1),
or under some natural restrictions (such as $\W$ being a convex set in Euclidean space).

The results of this paper concern the direction ``spectral $\Rightarrow$
tile'' of the  Spectral Set Conjecture in finite Abelian groups and $\Z^d$
and
$\R^d$. In Section \ref{sec2} we prove the validity of this direction of the
conjecture
for sets of at most 5 elements (and the corresponding union of unit cubes in the case of $\R^d$).
On the other hand, in Section \ref{sec3} we
show a particular
example of a
6-element set in $G=\Z_8^3$ which is spectral but does not tile $G$ by
translation. This allows us to improve  the results of \cite{tao} and
\cite{mat}, by producing a counterexample in $\R^3$.
Finally in Section \ref{sec:algorithmic} we discuss the algorithmic aspect
of when a given set tiles a finite cyclic group or is spectral in that group.
Such questions arise naturally and, unfortunately, with no satisfying answer,
when one searches for counterexamples to the Spectral Set Conjecture
in cyclic groups (any such counterexample would extend to a counterexample in $\R$).

\section{Spectral sets of small size}\label{sec2}

We recall the following notations and definitions:

A $k\times k$ complex matrix $H$ is called a  (complex) {\it Hadamard matrix} if all
entries of $H$ have
absolute value 1, and $HH^\ast =kI$ (where $I$ denotes the identity
matrix). This means that the rows (and also the columns) of $H$ form an
orthogonal basis of $\Co^k$.  A {\it log-Hadamard matrix} is any real square matrix
$(h_{i,j})_{i,j=1}^k$ such that the matrix $(e^{2\pi i h_{i,j}})_{i,j=1}^k$ is
Hadamard.

Based on the ideas of \cite{tao} and \cite{mat} we will make use of the
connection between complex Hadamard matrices and spectral sets to prove the following:

\begin{theorem}\label{finite}
Let $G$ be any finite Abelian group and $A\subset G$ a spectral set in $G$ with $|A|\le 5$. Then
$A$ tiles $G$.
\end{theorem}
\begin{proof}
We begin the proof by recalling a general property of spectral sets and tiles. Namely,
if $H\le G$ is a subgroup and $A\subset H$, then $A$ is spectral (resp. a tile) in $G$ if and
only if $A$ is spectral (resp. a tile) in $H$.
The only non-trivial statement here is that if
$A$ is spectral in $H$ then it is also spectral in $G$. This follows directly from  
the fact that each
character of $H$ can be extended, not necessarily uniquely, to a character of $G$. To see this, 
we need to show that the natural (restriction) homomorphism  $\phi : \hat{G} \to \hat{H}$ is 
surjective. Any element in the kernel of $\phi$ is constant on each coset of $H$, therefore
it can be identified with a character of the factor group $G/H$. Then, the surjectivity of $\phi$ 
follows from comparing the number of elements of the groups $G, H, G/H$.

It is a well known fact that for any finite Abelian group $G$ we may choose natural numbers $N, d$ such that
$G\le \Z_N^d$.
By the remark above, it is enough to prove the statement for groups of the type $G=\Z_N^d$.
This observation makes the proof technically simpler.

The essential part of the proof relies on the fact that we have a full characterization of
complex Hadamard matrices up to order 5.

We identify the elements of $G$ and $\hat{G}$, the group of characters of $G$,
with $d$-dimensional column- and row-vectors, respectively.
Let $A\subset G=\Z_N^d$, with $|A|=k\le 5$. We regard $A$ as a $d\times k$ matrix
with integer coefficients.
If $L\subset \hat{G}$ is a spectrum of $A$ (regarded as a
$d\times k$ matrix), then
$H:=\frac{1}{N}L\cdot A$ is log-Hadamard (where the matrix multiplication can be taken $\bmod \ N$).
Multiplication by the matrix $L$ defines a homomorphism from $G$ to $\Z_N^k$, and the images of
the elements of $A$ are given by the columns $c_j$ of $L\cdot A$ ($0\le j\le k-1$).  By
Proposition
2.3 in \cite{mat}
it is enough to prove that the vectors $c_j$ tile $\Z_N^k$.

In the cases  $k=1,2,3,5$ this follows immediately from the uniqueness (up to natural
equivalence) of complex
Hadamard matrices of order $k$.
This uniquness is trivial for $k=1,2,3$, while the case $k=5$
is settled in \cite{haag}.
Indeed, we can assume without loss of
generality that $0\in A$ and $0\in L$ (due to the trivial translation invariance of the notion of
spectrality and spectrum), and this already implies that the matrix $H=h_{m,j}$
is (after a permutation of columns) given by $h_{m,j}=\frac{1}{k}mj$ ($0\le m,j\le k-1$).
It follows that $N$ is a multiple of $k$, say $N=Mk$, and the column vectors $c_j$ are given as
$c_j=(0, jM, 2jM, \dots ,(k-1)jM)^T$. In order to see that these vectors tile $\Z_N^d$ we invoke
Proposition 2.3 in \cite{mat} once again. Let $V:=(0,1,0,\dots ,0)$, and consider the $\bmod \ N$
product $V\cdot L\cdot A=(0,M,2M, \dots (k-1)M)$. It is obvious that the set
$\{0,M,2M, \dots (k-1)M\}$ tiles $\Z_N$, and therefore the columns $c_j$ tile $\Z_N^k$.

The case $k=4$ is settled in a similar manner, although we have no uniqueness of Hadamard matrices
in this case. The general form of a $4\times 4$ complex Hadamard matrix is given (see e.g. 
\cite{haag}, Proposition 2.1) by the
parametrization
$$
U=\left (
\begin{array}{cccc}
1&1&1&1\\
1&1&-1&-1\\
1&-1&e^{2\pi i\phi}&-e^{2\pi i\phi}\\
1&-1&-e^{2\pi i\phi}&e^{2\pi i\phi}
\end{array}
\right ).
$$

Due to the presence of -1's it follows that $N$ must be a multiple of 2, say $N=2M$. The matrix
is $LA$ is then given by
$$
LA=\left (
\begin{array}{cccc}
0&0&0&0\\
0&0&M&M\\
0&M&M+t&t\\
0&M&t&M+t
\end{array}
\right ).
$$

To see that the columns of this matrix tile $\Z_N^4$, consider the matrix
$$
V_2:=\left (
\begin{array}{cccc}
0&1&0&0\\
0&0&1&0
\end{array}
\right ).
$$
Then
$$
V_2LA=\left (
\begin{array}{cccc}
0&0&M&M\\
0&M&M+t&t
\end{array}
\right ).
$$

It is trivial to check that the columns of this matrix tile $\Z_N^2$, and, by Proposition 2.3 in
\cite{mat}, this implies that the columns of $LA$ tile $\Z_N^4$.

\end{proof}

Next, we extend the previous result to the infinite grid $\Z^d$. First, we need to establish the
rationality of the spectra in the cases considered.

\begin{proposition}\label{rational}
Let $A\subset Z^d$ be a spectral set with $|A|\le 5$. Then $A$ admits a rational spectrum.
\end{proposition}
\begin{proof}
Note that we do not claim that {\it all} spectra of $A$ must be rational, but only that the
spectrum can be chosen rational.

The proof is an easy argument from linear algebra. Let us first consider the case $|A|=5$
(the cases $|A|=1,2,3$ are settled the same way, while $|A|=4$ will require some extra
considerations). Let $L\subset \T^d$ denote a spectrum of $A$.
We may assume that $0\in A$ and
$0\in L$.
Then, after a permutation of elements of $A$, we have
$$
LA=\frac{1}{5}\left (
\begin{array}{ccccc}
0&0&0&0&0\\
0&1&2&3&4\\
0&2&4&1&3\\
0&3&1&4&2\\
0&4&3&2&1
\end{array}
\right )   \ \ (\bmod \ 1).
$$
Let $l_{m,j}$ denote the elements of $L$.
Considering, for example, the second row of $L$ we see that
there exist integers $z_{2,1}, \dots ,z_{2,5}$ such that
$$
(l_{2,1},\dots ,l_{2,d})\cdot A= \left(z_{2,1}, z_{2,2}+\frac{1}{5}, \ldots z_{2,5}+\frac{4}{5}\right).
$$
Regarding this equation as a set of linear equations with variables $l_{2,1},\ldots ,l_{2,d }$ we
see that if there exists a solution, then the solution can be chosen rational. The same argument
holds for the other rows of $L$.

We now turn to the case $|A|=4$. Then, for some $q\in [0,1]$ we have
$$
LA=\left (
\begin{array}{cccc}
0&0&0&0\\
0&0&\frac{1}{2}&\frac{1}{2}\\
0&\frac{1}{2}&\frac{1}{2}+q&q\\
0&\frac{1}{2}&q&\frac{1}{2}+q
\end{array}
\right ) (\bmod \ 1).
$$
If $q$ is rational then the previous argument applies. If $q$ is irrational we need some
additional considerations. Applying the previous argument we see that
the first two rows of $L$ can be chosen rational even in this case (they do not depend on $q$).
It is also clear that the fourth row of $L$ can be chosen as the sum of the second and third
rows. Consider therefore the third row only. For some integers $z_{3,1}, \ldots ,z_{3,4}$ we have
$$
(l_{3,1},\ldots ,l_{3,d})\cdot A= \left(z_{3,1}, z_{3,2}+\frac{1}{2}, z_{3,3}+\frac{1}{2}+q, z_{3,4}+q\right)
$$
Regarding this equation as a set of linear equations with variables $l_{2,1},\ldots ,l_{2,d}, q$
we see that if there exists a solution then all the variables (including $q$)
can be chosen rational.
\end{proof}

Now, we are in position to prove the analogue of Theorem \ref{finite} in $\Z^d$.

\begin{theorem}\label{integers}
Let $A\subset \Z^d$ be a spectral set in $\Z^d$ with $|A|\le 5$. Then
$A$ tiles $\Z^d$.

\end{theorem}
\begin{proof}
By Proposition \ref{rational} we can choose the spectrum of $A$ rational. This means that $A$ is
already spectral in some finite group $G=\Z_{n_1}\times\cdots\times\Z_{n_d}$.
By Theorem \ref{finite} we conclude that $A$ tiles
$G$, and therefore it tiles $\Z^d$.
\end{proof}

Finally, we formulate the corresponding statement in $\R^d$.

\begin{theorem}\label{reals}
Let $A$ be a set in $\Z^d$ with $|A|\le 5$. Let $\W:=A+[0,1)^d$ be the corresponding union of unit
cubes. If $\W$ is spectral then it tiles $\R^d$.
\end{theorem}
\begin{proof}
By Theorem 4.2 in \cite{kolmat} the set $\W$ is spectral in $\R^d$ if and only if $A$ is spectral
in $\Z^d$. Therefore, Theorem \ref{integers} implies the result.
\end{proof}

\section{The Spectral Set Conjecture is false in dimension 3} \label{sec3}

In the previous section we saw how a characterization of Hadamard
matrices (up to order 5)
leads to an understanding of tiling properties of spectral sets.
The situation changes drastically when we consider Hadamard matrices of order
6. On the one hand side we do not seem to have a full characterizaton
of complex Hadamard matrices of order 6 or greater. On the
other hand we do possess particular examples (and even a description
\cite{dita} of some
parametric families of $6\times 6$ Hadamard matrices) which show that the
statement of Theorem \ref{finite} does not hold for sets of size $|A|=6$.
Such an example was exploited in \cite{tao} (and subsequently in \cite{mat})
to disprove
Fuglede's conjecture in $\Z_3^5$ (resp. $\Z_3^4$) and then transfer the
counterexample to $\Z^5$ and $\R^5$ (resp. $\Z^4$ and $\R^4$).

In this section we improve the results mentioned above by showing that the
direction ``spectral $\Rightarrow$ tile'' of
the Spectral Set Conjecture fails already in dimension 3.
We use a combination of the ideas of \cite{tao} and \cite{mat}.

Our starting point is the following numerical example of a particular $6\times
6$ log-Hadamard matrix.
Let
$$
H:=\frac{1}{8}\left (
\begin{array}{cccccc}
0&0&0&0&0&0\\
0&4&2&6&6&2\\
0&2&4&1&5&6\\
0&6&3&4&2&7\\
0&6&7&2&4&3\\
0&2&6&5&1&4
\end{array}
\right ).
$$
It is easy to check that $H$ is log-Hadamard (this example corresponds to a
particular
element of the one-parameter family of $6\times 6$ Hadamard matrices given in
\cite[p.\ 5357]{dita}).

Next, we observe that the ``mod 8 rank'' of $8H$ is 3. This means that there
exist integer matrices $A$ and $L$ of size $3\times 6$ and $6\times 3$, respectively,
such that $8H=LA$ ($\bmod \ 8$). A possible example of such a
decomposition is the following:
$$
A:=\left (
\begin{array}{cccccc}
0&2&4&1&5&6\\
0&6&3&4&2&7\\
0&6&7&2&4&3
\end{array}
\right ) \ \ \ \ \mathrm{and}  \ \ \
L:=\left (
\begin{array}{cccccc}
0&0&0\\
0&1&1\\
1&0&0\\
0&1&0\\
0&0&1\\
7&1&1
\end{array}
\right ).
$$
   This means that the set (of the columns of) $A$ is spectral in the group
$G=\Z_8^3$. However, $|A|=6$, therefore $A$ cannot tile  $G$ due to obvious
divisibility reasons.

We can transfer this counterexample to $\Z^3$ and $\R^3$ in the same way as
described in \cite{tao} and \cite{mat}. Namely, the set $A_n:=A+8[0,n)^3$ is spectral in
$\Z^3$ for every $n$, but it does not tile $\Z^3$ for large enough values of
$n$ (see Propositions 2.1 and 2.5 in \cite{mat}).
Finally, for such a large value $n$, the set $\W=A_n+[0,1)^3$ is a finite
union of unit cubes, which is spectral in $\R^3$ but does not tile $\R^3$ (cf.\
the proof of Theorem 3.1 in \cite{mat}).

We see that the counterexample is simply based on the existence of a
particular log-Hadamard
matrix with prescribed properties (namely, to get the non-divisibility
condition the size (6) must not divide the
common denominator of the entries (8); the dimension (3) is then determined by the
'rank' of the matrix, i.e. the smallest possible decomposition).
In principle, this would seem to be the way to search for counterexamples in
dimensions 2 and 1.

However, due to the lack of characterization of Hadamard matrices of order
greater than 5, we do not know of any methodical way of {\it how} to produce
Hadamard
matrices with the desired properties. The particular example above was simply
found in the literature by looking through the existing (fairly scattered)
results on complex Hadamard matrices.

It seems that in order to prove or disprove this direction of the Spectral
Set
Conjecture in
dimensions 1 and 2, some new ideas will be needed. Let us also mention that
the direction ``tile $\Rightarrow$ spectral'' has recently been disproved
\cite{kolmat} in dimension 5 and higher, but remains open in lower dimensions.

\section{Algorithmic questions}\label{sec:algorithmic}
In searching for counterexamples to the Spectral Set Conjecture using a computer one is
immediately faced with the following problem: it appears that no efficient
algorithm is known to decide if a given subset $E$ of some group, say the cyclic group $\Z_n$, is a tile.
The customary definition of efficiency for such an algorithm demands that it should run in time polynomial in $n$.
No efficent algorithm appears to be known also in the case of deciding if a given set
is spectral or not.
The obvious algorithms, namely to test all possible tiling complements and all possible spectra,
are clearly exponential in $n$.

Both questions would easily be solvable given an efficient algorithm
to decide the following question.

\begin{quote}
{\sc Problem DIFF}:
\begin{itemize}
\item INPUT: We are given a positive integer $n$ and sets $E, D \subseteq \Z_n$
\item OUTPUT: Maximal $k$ such that there exists $A \subseteq E$, with $\Abs{A} = k$,
and $A-A \subseteq D$.
\end{itemize}
\end{quote}

Indeed, assuming one could solve the above problem in polynomial time in $n$,
one could use this postulated procedure to decide if $T\subseteq\Z_n$ tiles $\Z_n$ by using
$E = \Z_n$ and $D = (T-T)^c \cup \Set{0}$.
If $x = n / \Abs{T}$ is an integer and the answer to this query is equal to $x$ then $T$ tiles $\Z_n$ with translation set
equal to the set $A$ in the description of problem DIFF, otherwise $T$ is not a tile.

To decide spectrality of $S \subseteq \Z_n$ in polynomial time, first compute
the Fourier Transform of $\chi_S$, its zero set $Z$ and invoke the procedure
for problem DIFF with $E = D = \Set{0} \cup Z$.
If the answer $k$ is equal to $\Abs{S}$ then, and only then, $S$ is spectral in $\Z_n$
and one of its spectra is the set $A$ in the description of problem DIFF.

It is easily seen that {\sc DIFF} is solvable in polynomial time if and only if the corresponding
decision problem, described below, is solvable in polynomial time.

\begin{quote}
{\sc Problem DIFF'}:
\begin{itemize}
\item INPUT: We are given positive integers $k$ and $n$ and sets $E, D \subseteq \Z_n$
\item OUTPUT: YES if there exists $A \subseteq E$, with $\Abs{A} = k$,
and $A-A \subseteq D$, and NO otherwise.
\end{itemize}
\end{quote}

We show now that problem {\sc DIFF'} is NP-complete \cite{garey-johnson}.
We remark however that this does not prove that the corresponding
tiling and spectrality decision problems are NP-complete.
\begin{theorem}\label{th:np-complete}
The decision problem {\sc DIFF'} is NP-complete.
\end{theorem}
\begin{proof}
The following problem (decision version of the maximum independent set optimization problem)
is well known \cite{garey-johnson} to be NP-complete:
\begin{quote}
{\sc Problem IND}:
\begin{itemize}
\item INPUT: We are given a simple graph on $n$ vertices and a positive integer $k$.
\item OUTPUT: YES if there exists a set of vertices $A$ none of which is connected
via an edge to another vertex in $A$ (the set $A$ is then called an independent set) and such that $\Abs{A}= k$.
\end{itemize}
\end{quote}
To prove that {\sc DIFF'} is NP-complete it suffices to give a polynomial time
algorithm to solve the arbitrary instance of {\sc IND} using a black-box algorithm that
solves {\sc DIFF'} in unit time.

Let $V = \Set{1,\ldots, n}$ be the set of vertices of a given graph $G$.
The first step is to define an efficiently computable function $\phi: V \to \Z_m$,
with $m$ bounded by a polynomial in $n$,
and a set $A \subseteq \Z_m$ such that
\beql{graph-map}
i,j \in V\ \mbox{are connected in $G$} \Leftrightarrow
 \phi(i)-\phi(j) \in A.
\end{equation}
For this we note that one can easily construct
such a map $\phi$ and set $A$ using the greedy method, as long as $m$
is large compared to $n$.
Indeed, choose $\phi(1) = 0$, and, having defined
$\phi(1),\ldots, \phi(r)$, $r<n$, define
$\phi(r+1)$ to be the minimum $\nu \in \Set{1,\ldots,m-1}$
such that
$$
\nu \notin \Set{\phi(i)+\phi(j)-\phi(l):\ 1\le i,j,l \le r}.
$$
(It is easy to see that this procedure can be carried out if $m$ is large enough, namely $m \ge C n^3$.)
This ensures that all differences $\phi(i)-\phi(j)$,
$1\le i\neq j \le n$ are distinct.

Next we define
$$
A = \Set{\phi(i)-\phi(j):\ \mbox{$i$ and $j$ are connected in $G$}}.
$$
To decide if $G$ has an independent set of size $k$ it is enough to use the
algorithm for {\sc DIFF'} for $\Z_m$ with the set $E = \phi(\Set{1,\ldots,n})$ and
$A$ as constructed above.
\end{proof}

\end{document}